\documentclass[12pt]{article}
\usepackage[reqno]{amsmath}
\usepackage{amssymb, theorem, enumerate}
\usepackage{graphicx}
\usepackage{array}

\usepackage[svgnames]{xcolor}
\usepackage{tikz}
\usetikzlibrary{decorations.markings}
\usetikzlibrary{shapes.geometric}

\pgfdeclarelayer{edgelayer}
\pgfdeclarelayer{nodelayer}
\pgfsetlayers{edgelayer,nodelayer,main}

\tikzstyle{none}=[inner sep=0pt]

\tikzstyle{rn}=[circle,fill=Red,draw=Black,line width=0.8 pt]
\tikzstyle{gn}=[circle,fill=Lime,draw=Black,line width=0.8 pt]
\tikzstyle{yn}=[circle,fill=Yellow,draw=Black,line width=0.8 pt]
\tikzstyle{blackcirc}=[circle,fill=Black,draw=Black, scale= 0.5]
\tikzstyle{whitecirc}=[circle,fill=White,draw=Black, scale=0.5]
\tikzstyle{newstyle}=[circle,fill=White,draw=Black]

\tikzstyle{simple}=[-,draw=Black,line width=1.000]
\tikzstyle{arrow}=[-,draw=Black,postaction={decorate},decoration={markings,mark=at position .5 with {\arrow{>}}},line width=2.000]
\tikzstyle{tick}=[-,draw=Black,postaction={decorate},decoration={markings,mark=at position .5 with {\draw (0,-0.1) -- (0,0.1);}},line width=2.000]

\expandafter\ifx\csname pdfpagebox\endcsname\relax\else
\fi

\theoremstyle{change}
{\theorembodyfont{\slshape}
\newtheorem{theorem}{Theorem.}[section]

}
\theorembodyfont{\rmfamily}

\def\sqr#1#2{{\vbox{\hrule height.#2pt
    \hbox{\vrule width.#2pt height#1pt \kern#1pt
        \vrule width.#2pt}\hrule height.#2pt}}}
\def\eqed{\sqr53}
\def\qed{%
    \ifmmode\eqno\eqed
    \else\nobreak\ \hfill\eqed\medbreak\fi}

\newcommand\ints{{\mathbb Z}}

\newcommand\pmat[1]{\begin{pmatrix} #1 \end{pmatrix}}

\title{Large regular bipartite graphs with median eigenvalue $1$}

\author{Krystal Guo\thanks{Supported in part by NSERC PGS.}\\[1mm]
   {Department of Mathematics}\\{Simon Fraser University}\\{Burnaby, B.C. V5A 1S6} \\ \texttt{krystalg@sfu.ca}
    \and
   Bojan Mohar\thanks{Supported in part by an NSERC Discovery Grant (Canada), by the Canada Research Chair program, and by the Research Grant of ARRS (Slovenia).}~~\thanks{On leave from: IMFM \& FMF, Department of Mathematics, University of Ljubljana, Ljubljana, Slovenia.} \\[1mm]
 {Department of Mathematics}\\{Simon Fraser University}\\{Burnaby, B.C. V5A 1S6}\\ \texttt{mohar@sfu.ca}
 }

\begin{document}
\maketitle

\begin{abstract}
A recent result of one of the authors says that every connected subcubic bipartite graph that is not isomorphic to the Heawood graph has at least one, and in fact a positive proportion of its eigenvalues in the interval $[-1,1]$. We construct an infinite family of connected cubic bipartite graphs which have no eigenvalues in the open interval $(-1,1)$, thus showing that the interval $[-1,1]$ cannot be replaced by any smaller symmetric subinterval even when allowing any finite number of exceptions. Similar examples with vertices of larger degrees are considered and it is also shown that their eigenvalue distribution has somewhat unusual properties. By taking limits of these graphs, we obtain examples of infinite vertex-transitive $r$-regular graphs for every $r\ge3$, whose spectrum consists of points $\pm1$ together with intervals $[r-2,r]$ and $[-r,-r+2]$. These examples shed some light onto a question communicated by Daniel Lenz and Matthias Keller with motivation in relation to the Baum-Connes conjecture.

\vspace{5pt}

\vspace{5pt}
\end{abstract}

\section{Introduction}

Let $G$ be a graph of order $n$ and let $\lambda_1(G) \ge \lambda_2(G) \ge\cdots\ge\lambda_n(G)$ be the eigenvalues of its adjacency matrix. In this paper we are interested in the \emph{median eigenvalues} $\lambda_{\lfloor\frac{n+1}{2}\rfloor}(G)$ and $\lambda_{\lceil\frac{n+1}{2}\rceil}(G)$.

The eigenvalues of graphs can be useful descriptors of certain combinatorial properties of the graph. The most important one is related to the \emph{spectral gap} (the difference $\lambda_1(G) -\lambda_2(G)$), which certifies expansion properties of graphs, see, e.g. \cite{AM85, MP91, HLW06}. A more recent application comes from mathematical chemistry \cite{FP10a, FP10b, M13} in relation to the HOMO-LUMO energies of molecular graphs. In this setting, the median eigenvalues play the central role. Motivated by questions in mathematical chemistry, one of the authors proved the following curious result.

\begin{theorem}[Mohar \cite{MoharCPC}]
\label{thm:cubic bipartite}
Let $G$ be a bipartite subcubic graph. If every connected component of $G$ is isomorphic to the Heawood graph, then its median eigenvalues are $\pm\sqrt{2}$. In any other case, the median eigenvalues lie in the interval $[-1,1]$.
\end{theorem}

In fact, this theorem can be strengthened.

\begin{theorem}[Mohar \cite{MoharCPC}]
\label{thm:cubic bipartite linearly many}
There is a constant $\delta>0$ such that for every bipartite subcubic graph $G$ of order $n$, none of whose connected components is isomorphic to the Heawood graph, at least $\lceil\delta n\rceil$ of its eigenvalues belong to the interval $[-1,1]$.
\end{theorem}

However, the paper \cite{MoharCPC} leaves an open question:
\begin{quote}
\emph{Is there a strengthening of Theorem \ref{thm:cubic bipartite} where the interval $[-1,1]$ is replaced by a smaller symmetric interval around 0 if we allow a finite number of exceptional graphs?}
\end{quote}

In this note we answer the question in the negative by proving:

\begin{theorem}
\label{thm:median 1}
There are infinitely many connected cubic bipartite graphs that have no eigenvalues in the open interval $(-1,1)$.
\end{theorem}

By Theorem \ref{thm:cubic bipartite linearly many}, large graphs in the family of Theorem \ref{thm:median 1} will have $\pm 1$ as eigenvalues of large multiplicity.

The construction of graphs used to prove Theorem \ref{thm:median 1} can be generalized to larger vertex degrees, providing examples with unusual eigenvalue distributions.

\begin{theorem}
\label{thm:median large degree}
For every integer $r\ge3$, there are infinitely many connected $r$-regular bipartite graphs with median eigenvalues $\pm 1$ but with no eigenvalues in the intervals $(-1,1)$ and $\pm(1,r-2)$.
\end{theorem}

By taking limits of graphs from Theorem \ref{thm:median large degree}, we obtain examples of infinite vertex-transitive $r$-regular graphs for every $r\ge3$, whose spectrum consists of $\pm1$ together with intervals $[r-2,r]$ and $[-r,-r+2]$. These provide examples of infinite vertex-transitive graphs whose spectrum is not formed by a single interval. Daniel Lenz and Matthias Keller (private communication) studied the spectrum of regular tessellations of the hyperbolic plane. They conjectured that the spectrum consists of finitely many intervals with absolutely continuous spectrum and some eigenvalues, as this is the case for periodic operators on the plane. A partial answer for questions of this type come from discrete cases of the Baum-Connes conjecture (see \cite{Valette02,HLS02}), which gives examples whose spectrum consists of only one interval.

\section{Graphs $W_{n,k}$ and their eigenvalues}

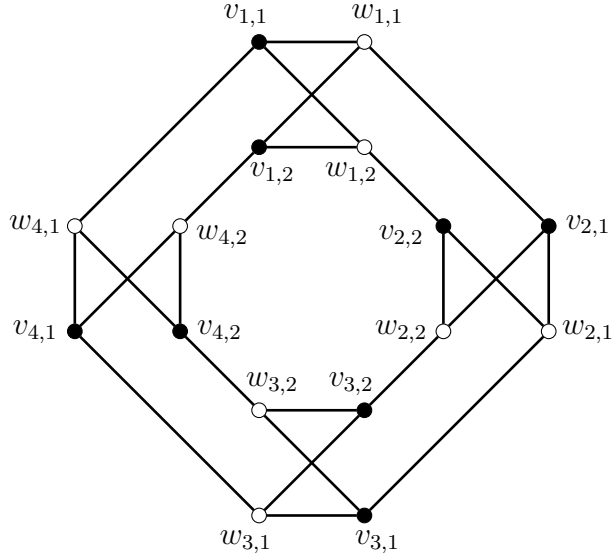
\begin{figure}[h]
\begin{center}
\begin{tikzpicture}[scale=0.7]
	\begin{pgfonlayer}{nodelayer}
		\node [style=blackcirc] (0) at (1, -4.5) {};
		\node [style=blackcirc] (1) at (-1, 4.5) {};
		\node [style=blackcirc] (2) at (2.5, 1) {};
		\node [style=blackcirc] (3) at (-1, 2.5) {};
		\node [style=blackcirc] (4) at (-2.5, -1) {};
		\node [style=blackcirc] (5) at (1, -2.5) {};
		\node [style=blackcirc] (6) at (-4.5, -1) {};
		\node [style=blackcirc] (7) at (4.5, 1) {};
		\node [style=whitecirc] (8) at (1, 4.5) {};
		\node [style=whitecirc] (9) at (1, 2.5) {};
		\node [style=whitecirc] (10) at (2.5, -1) {};
		\node [style=whitecirc] (11) at (4.5, -1) {};
		\node [style=whitecirc] (12) at (-1, -2.5) {};
		\node [style=whitecirc] (13) at (-1, -4.5) {};
		\node [style=whitecirc] (14) at (-2.5, 1) {};
		\node [style=whitecirc] (15) at (-4.5, 1) {};
		\node [style=none] (16) at (1.25, 5) {$w_{1,1}$};
		\node [style=none] (17) at (-1.25, 5) {$v_{1,1}$};
		\node [style=none] (18) at (-5.25, 1) {$w_{4,1}$};
		\node [style=none] (19) at (-5.25, -1) {$v_{4,1}$};
		\node [style=none] (20) at (-1.25, -5) {$w_{3,1}$};
		\node [style=none] (21) at (1.25, -5) {$v_{3,1}$};
		\node [style=none] (22) at (5.25, -1) {$w_{2,1}$};
		\node [style=none] (23) at (5.25, 1) {$v_{2,1}$};
		\node [style=none] (24) at (-0.75, 2) {$v_{1,2}$};
		\node [style=none] (25) at (-1.7, 0.8) {$w_{4,2}$};
		\node [style=none] (26) at (-1.75, -1) {$v_{4,2}$};
		\node [style=none] (27) at (-0.75, -2) {$w_{3,2}$};
		\node [style=none] (28) at (0.75, -2) {$v_{3,2}$};
		\node [style=none] (29) at (1.7, -1) {$w_{2,2}$};
		\node [style=none] (30) at (1.7, 0.8) {$v_{2,2}$};
		\node [style=none] (31) at (0.75, 2) {$w_{1,2}$};
	\end{pgfonlayer}
	\begin{pgfonlayer}{edgelayer}
		\draw [style=simple] (1) to (8);
		\draw [style=simple] (8) to (3);
		\draw [style=simple] (3) to (9);
		\draw [style=simple] (9) to (1);
		\draw [style=simple] (14) to (6);
		\draw [style=simple] (4) to (15);
		\draw [style=simple] (7) to (10);
		\draw [style=simple] (11) to (2);
		\draw [style=simple] (12) to (5);
		\draw [style=simple] (5) to (13);
		\draw [style=simple] (13) to (0);
		\draw [style=simple] (0) to (12);
		\draw [style=simple] (12) to (4);
		\draw [style=simple] (14) to (3);
		\draw [style=simple] (9) to (2);
		\draw [style=simple] (10) to (5);
		\draw [style=simple] (13) to (6);
		\draw [style=simple] (15) to (1);
		\draw [style=simple] (8) to (7);
		\draw [style=simple] (7) to (11);
		\draw [style=simple] (11) to (0);
		\draw [style=simple] (6) to (15);
		\draw [style=simple] (4) to (14);
		\draw [style=simple] (2) to (10);
	\end{pgfonlayer}
\end{tikzpicture}
\caption{The graph $W_{4,2}$ with bipartite classes shown as black and white vertices. \label{fig:G4}}
\end{center}
\end{figure}

We consider a family of $(k+1)$-regular graphs $W_{n,k}$ ($n\geq 2$, $k\ge 2$) on $2nk$ vertices that are defined as follows. The graph $W_{n,k}$ has vertices 
$ v_{i,j}$ $(1\le i\le n, 1\le j\le k)$
and 
$ w_{i,j}$ $ (1\le i\le n, 1\le j\le k)$.
Its edge-set consists of edges joining every $v_{i,j}$ and every $w_{i,l}$ ($1\le i\le n$, $1\le j\le k$, $1\le l\le k$) together with a matching consisting of all edges joining $w_{i,j}$ with $v_{i+1,j}$ ($1\le i\le n$, $1\le j\le k$, where the index $i+1$ is taken modulo $n$). The graph $W_{n,k}$ is clearly $(k+1)$-regular, bipartite and vertex-transitive. The graph $W_{4,2}$ is depicted in Figure \ref{fig:G4}.

\begin{theorem}\label{thm:evals}
For every $n \geq 2$ and $k \geq 2$, the spectrum of the graph $W_{n,k}$ consists of eigenvalues $\pm 1$, each with multiplicity $(k -1) n$, and the remaining eigenvalues are equal to $\pm \tau_j$ ($j=0, \ldots n-1$) each with multiplicity $2$, where
\[ \tau_j = \sqrt{k^2 +1 + 2k\cos\left(\tfrac{2\pi j}{n}\right)} .\]
\end{theorem}

The stated multiplicity of $\pm 1$ may be higher if $\tau_j =1$ for some $j$, which occurs when $k=2$ and $j = \frac{n}{2}$. Whenever $j + \ell = n$, we have that $\tau_j = \tau_{\ell}$, so that the eigenvalues equal to $\tau_j$ occur with multiplicity $4$. 

The graphs $W_{n,k}$ have their median eigenvalues equal to $\pm 1$. To see this, let $z = \tau_j$. Then 
\[
z^2 =  k^2 + 1  + 2k\cos\left(\tfrac{2\pi j}{n}\right).
\]
Since
$ -1 \leq \cos\left(\frac{2\pi j}{n}\right) \leq 1 $
we have that 
\[
k^2 + 1 - 2k \leq z^2 \leq k^2 + 1 + 2k,
\]
which gives that
\[
k-1 \leq z \leq k+1.
\] 
This shows that all eigenvalues of $W_{n,k}$ which are different from $\pm 1$ lie in the intervals $[k-1, k +1]$ and $[-k-1, -k+1]$, so they are strangely concentrated when $k$ is large. This, in particular, proves Theorem \ref{thm:median 1} (by taking $k=2$) and Theorem \ref{thm:median large degree} (by taking $r = k+1$). The rest of this section is devoted to the proof of Theorem \ref{thm:evals}. 

\vspace{10pt}
\noindent \textsl{Proof of Theorem \ref{thm:evals}: }  Let $A = A(W_{n,k})$ be the adjacency matrix of $W_{n,k}$. We will find the eigenvalues of $A$ by considering those of $A^2$. Since $W_{n,k}$ is bipartite, its eigenvalues are symmetric about the origin. Then, all eigenvalues of $A^2$ are non-negative and every eigenvalue $\mu$ of $A^2$ corresponds to two eigenvalues $\sqrt{\mu}$ and $-\sqrt{\mu}$ of $A$, whose multiplicities are the same as that of $\mu$. 

\begin{figure}[h]
\begin{center}
\begin{tikzpicture}[scale=0.62]
	\begin{pgfonlayer}{nodelayer}
		\node [style=blackcirc] (0) at (-5.5, 4) {};
		\node [style=blackcirc] (1) at (-5.5, 2) {};
		\node [style=blackcirc] (2) at (-5.5, -2) {};
		\node [style=blackcirc] (3) at (-5.5, -4) {};
		\node [style=blackcirc] (4) at (-3.5, -0) {};
		\node [style=blackcirc] (5) at (-1.5, -0) {};
		\node [style=blackcirc] (6) at (-7.5, -0) {};
		\node [style=blackcirc] (7) at (-9.5, -0) {};
		\node [style=whitecirc] (8) at (5.5, -4) {};
		\node [style=whitecirc] (9) at (3.5, -0) {};
		\node [style=whitecirc] (10) at (9.5, -0) {};
		\node [style=whitecirc] (11) at (5.5, 2) {};
		\node [style=whitecirc] (12) at (7.5, -0) {};
		\node [style=whitecirc] (13) at (5.5, 4) {};
		\node [style=whitecirc] (14) at (1.5, -0) {};
		\node [style=whitecirc] (15) at (5.5, -2) {};
		\node [style=none] (16) at (-5.5, 4.75) {$v_{1,1}$};
		\node [style=none] (17) at (-5.5, -4.75) {$v_{3,1}$};
		\node [style=none] (18) at (-0.6, -0.1) {$v_{2,1}$};
		\node [style=none] (19) at (-10.5, -0) {$v_{4,1}$};
		\node [style=none] (20) at (-5.5, 1.25) {$v_{1,2}$};
		\node [style=none] (21) at (-4.3, -0) {$v_{2,2}$};
		\node [style=none] (22) at (-6.7, -0) {$v_{4,2}$};
		\node [style=none] (23) at (-5.5, -1.3) {$v_{3,2}$};
		\node [style=none] (24) at (5.5, 4.75) {$w_{1,1}$};
		\node [style=none] (25) at (0.6, 0.1) {$w_{4,1}$};
		\node [style=none] (26) at (5.5, -4.75) {$w_{3,1}$};
		\node [style=none] (27) at (10.5, -0) {$w_{2,1}$};
		\node [style=none] (28) at (5.5, 1.25) {$w_{1,2}$};
		\node [style=none] (29) at (5.5, -1.3) {$w_{3,2}$};
		\node [style=none] (30) at (6.6, -0) {$w_{2,2}$};
		\node [style=none] (31) at (4.4, -0) {$w_{4,2}$};
	\end{pgfonlayer}
	\begin{pgfonlayer}{edgelayer}
		\draw [bend left=45, looseness=1.00] (0) to (5);
		\draw [bend left=45, looseness=1.00] (5) to (3);
		\draw [bend left=45, looseness=1.00] (3) to (7);
		\draw [bend left=45, looseness=1.00] (7) to (0);
		\draw [bend right, looseness=1.00] (0) to (1);
		\draw [bend right, looseness=1.00] (4) to (5);
		\draw [bend right, looseness=1.00] (2) to (3);
		\draw [bend right, looseness=1.00] (6) to (7);
		\draw [bend left, looseness=1.00] (6) to (0);
		\draw [bend left, looseness=1.00] (0) to (4);
		\draw [bend left, looseness=1.00] (4) to (3);
		\draw [bend left, looseness=1.00] (3) to (6);
		\draw [bend left=45, looseness=0.75] (6) to (1);
		\draw [bend left, looseness=1.00] (1) to (5);
		\draw [bend left, looseness=1.00] (5) to (2);
		\draw [bend left, looseness=1.00] (2) to (7);
		\draw [bend right=45, looseness=0.75] (6) to (2);
		\draw [bend right=45, looseness=0.75] (2) to (4);
		\draw [bend right=45, looseness=0.75] (4) to (1);
		\draw [bend right, looseness=1.00] (1) to (7);
		\draw [bend left, looseness=1.00] (0) to (1);
		\draw [bend left, looseness=1.00] (4) to (5);
		\draw [bend left, looseness=1.00] (2) to (3);
		\draw [bend right, looseness=1.00] (7) to (6);
		\draw [bend left=45, looseness=1.00] (13) to (10);
		\draw [bend left=45, looseness=1.00] (10) to (8);
		\draw [bend left=45, looseness=1.00] (8) to (14);
		\draw [bend left=45, looseness=1.00] (14) to (13);
		\draw [bend right, looseness=1.00] (13) to (11);
		\draw [bend right, looseness=1.00] (12) to (10);
		\draw [bend right, looseness=1.00] (15) to (8);
		\draw [bend right, looseness=1.00] (9) to (14);
		\draw [bend left, looseness=1.00] (9) to (13);
		\draw [bend left, looseness=1.00] (13) to (12);
		\draw [bend left, looseness=1.00] (12) to (8);
		\draw [bend left, looseness=1.00] (8) to (9);
		\draw [bend left=45, looseness=0.75] (9) to (11);
		\draw [bend left, looseness=1.00] (11) to (10);
		\draw [bend left, looseness=1.00] (10) to (15);
		\draw [bend left, looseness=1.00] (15) to (14);
		\draw [bend right=45, looseness=0.75] (9) to (15);
		\draw [bend right=45, looseness=0.75] (15) to (12);
		\draw [bend right=45, looseness=0.75] (12) to (11);
		\draw [bend right, looseness=1.00] (11) to (14);
		\draw [bend left, looseness=1.00] (13) to (11);
		\draw [bend left, looseness=1.00] (12) to (10);
		\draw [bend left, looseness=1.00] (15) to (8);
		\draw [bend right, looseness=1.00] (14) to (9);
	\end{pgfonlayer}
\end{tikzpicture}
\caption{The graph $M$ corresponding to $W_{4,2}$. \label{fig:M4}}
\end{center}
\end{figure}
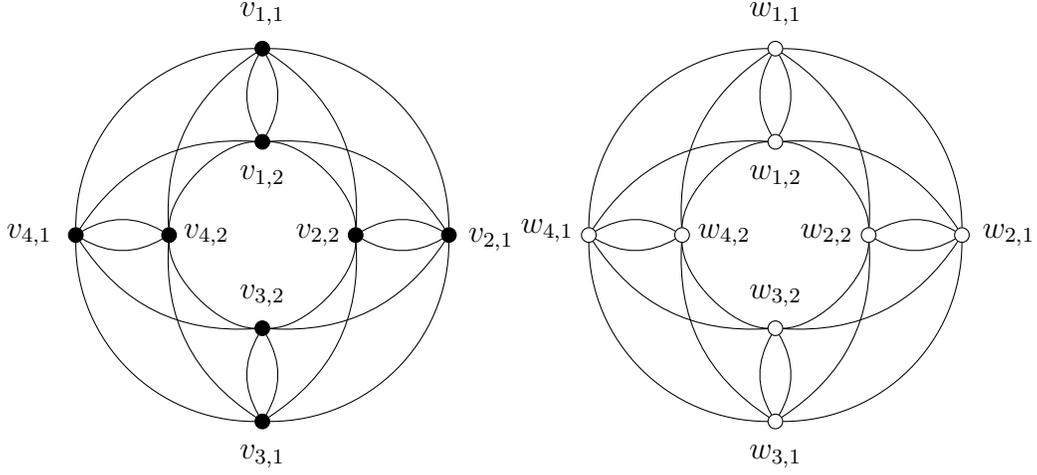

Since $W_{n,k}$ is $(k+1)$-regular, the diagonal entries of $A^2$ are all equal to $k+1$. We may then consider the loopless multigraph $M$ with adjacency matrix $A' = A^2 - (k+1)I_{2nk}$, where $I_{j}$ denotes the $j\times j$ identity matrix. We will use the fact that the $u,v$-entry of $A^2$ is equal to the number of walks of length two between vertices $u$ and $v$ in $W_{n,k}$. Since $W_{n,k}$ is bipartite, there are no walks of length two between vertices in different parts of the bipartition. Then $M$ consists of two connected components with vertex sets $V$ and $W$ where
\[
V = \{v_{i,j }\mid 1\le i\le n,\, 1\le j\le k \}
\]
and
\[
W = \{w_{i,j }\mid 1\le i\le n,\, 1\le j\le k \}.
\]
In $M$, for each $v_{i,j}$, there is one edge starting at $v_{i,j}$ and terminating at each of $v_{i+1, j'}$ and $v_{i-1, j'}$ for $1\le j' \le k$, where the indices are understood to be modulo $n$. There are $k$ edges with endpoints $v_{i,j}$ and $v_{i,j'}$ for $1\le j' \le k$ and $j \neq j$. In particular, each vertex has degree $2k + (k-1)k = k^2+k$ in $M$. The subgraphs of $M$ induced by $V$ and $W$ are isomorphic and $M$ is a $(k^2 + k)$-regular multigraph with two connected components. The multigraph $M$ corresponding to $W_{4,2}$ is depicted in Figure \ref{fig:M4}. 
 
Since $M$ has two components, we may write 
\[
A' = \pmat{A_V & 0 \\ 0 &A_W}
\]
where $A_V$ is the adjacency matrix of the subgraph of $M$ induced by $V$ and $A_W$ is the adjacency matrix of the subgraph of $M$ induced by $W$. Note that $A_V$ and $A_W$ are cospectral, so we will now focus on finding the eigenvalues of $A_V$. We may further write $A_V$ as 
\[ A_V = B_1 + kB_2
\]
where $B_1$ has the adjacencies between each $v_{i,j}$ and $v_{i-1, j'}$ and $v_{i+1, j'}$ and $B_2$ records the adjacencies between each $v_{i,j}$ and $v_{i, j'}$. 

Let 
\[
\alpha = \sqrt{\tfrac{k}{2} + \tfrac{1}{2}\sqrt{k^2 - 4}}
\quad \text{and}  \quad
\beta = \tfrac{1}{\alpha}.
\]
Next, we let $Q$ be the $n \times nk$ matrix given by
\[ \setcounter{MaxMatrixCols}{20}
Q 
= \pmat{\alpha & \alpha & \ldots & \alpha & 0 & 0 & \ldots & 0 & 0 & 0 & \ldots & 0 & \ldots & \beta & \beta & \ldots & \beta \\
\beta & \beta & \ldots & \beta & \alpha & \alpha & \ldots & \alpha & 0 & 0 & \ldots & 0 & \ldots  & 0 & 0 & \ldots & 0 \\
0 & 0 & \ldots & 0 &\beta & \beta & \ldots & \beta &  \alpha & \alpha & \ldots & \alpha & \ldots  & 0 & 0 & \ldots & 0  \\
0 & 0 & \ldots & 0 & 0 & 0 & \ldots & 0 & \beta & \beta & \ldots & \beta & \ldots  & 0 & 0 & \ldots & 0 \\
\vdots & \vdots &  & \vdots & \vdots & \vdots &  & \vdots & \vdots & \vdots &  & \vdots  & \ddots  & \vdots& \vdots &  & \vdots \\
0 & 0 & \ldots & 0 & 0 & 0 & \ldots & 0 & 0 & 0 & \ldots & 0 & \ldots  & \alpha & \alpha & \ldots & \alpha 
}
\]
It is easy to see that
\[
QQ^T = k(\alpha^2 + \beta^2)I_n + k\alpha\beta A(C_n),
\]
where $C_n$ is the $n$-cycle. On the other hand,
\[
Q^TQ = (\alpha^2 + \beta^2) I_{nk} + \alpha\beta B_1 + (\alpha^2 + \beta^2) B_2.
\]
Observe that $\alpha$ and $\beta$ have been carefully chosen such that
$
\alpha\beta = 1
$
and
\[
\begin{split}
\alpha^2 + \beta^2 &= \frac{k + \sqrt{k^2 - 4}}{2} + \frac{2}{k + \sqrt{k^2 - 4}}   \\
&= \frac{(k + \sqrt{k^2 - 4})^2 + 4}{2(k + \sqrt{k^2 - 4})} \\
&= k.
\end{split} 
\]
Thus, we have that 
\[
QQ^T = k^2 I_n + kA(C_n)
\]
and 
\[
Q^TQ = k I_{nk} + B_1 + k B_2 = k I_{nk} + A_V.
\]
We will use $\sigma(N)$ to denote the multiset of eigenvalues of matrix $N$, and in the lists of the eigenvalues, the multiplicities will in the superscripts if greater than one. We have that 
\[
\sigma(QQ^T) = \{ k^2 + 2k \cos({\textstyle \frac{2\pi j}{n}}) \mid j = 0 ,\ldots, n-1 \}
\]
and 
\[
\sigma(Q^TQ) = \sigma(QQ^T)\cup \{0^{(nk-n)}\}.
\]
Therefore,
\[
\sigma(A_V) = \{ k^2 + 2k \cos({\textstyle \frac{2\pi j}{n}}) - k \mid j = 0 ,\ldots, n-1 \}\cup \{-k^{(nk-n)}\}
\]
and thus
\[
\sigma(A') = \{ (k^2 + 2k \cos({\textstyle \frac{2\pi j}{n}}) - k)^{(2)} \mid j = 0 ,\ldots, n-1 \}\cup \{-k^{(2nk-2n)}\}.
\]
From this, we obtain that
\[ 
\sigma(A^2) = \{ (k^2 + 2k \cos({\textstyle \frac{2\pi j}{n}}) + 1)^{(2)} \mid j = 0 ,\ldots, n-1 \}\cup \{1^{(2nk-2n)}\}. 
\]
Finally,
\[
\sigma(A) = \Big\{ \pm \sqrt{k^2 + 1 + 2k \cos({\textstyle \frac{2\pi j}{n}} )} \mid j = 0 ,\ldots, n-1 \Big\}\cup \{1^{(nk-n)}\}\cup \{-1^{(nk-n)}\}
\]
as claimed. \qed

\section{Extension to infinite graphs}\label{sec:inf}

Concerning the spectrum of infinite graphs, we refer to \cite{MW89} or \cite{Woess00}. For every fixed $k$, the limit of the graphs $W_{n,k}$, when $n$ tends to infinity, gives rise to an infinite vertex-transitive graph $Z_k$, whose vertices are $\hat{v}_{i,j}$ and $\hat{w}_{i,j}$, where $ i\in \ints$ and $j \in \{i,\ldots, k\}$. The edge set of $Z_k$ consists of edges joining every $\hat{v}_{i,j}$ and every $\hat{w}_{i,l}$ ($i \in \ints$, $1\le j\le k$, $1\le l\le k$) together with a matching consisting of all edges joining $\hat{w}_{i,j}$ with $\hat{v}_{i+1,j}$ ($i\in \ints$, $1\le j\le k$). If we remove $k$ vertices $w_{n,1}, \ldots, w_{n,k}$ from $W_{n,k}$, we obtain a graph $P_{n,k}$ which is isomorphic to the ball of radius $n-1$ around a vertex in $Z_k$. Here we take the \emph{ball of radius $r$ around $v$} as the subgraph induced by all vertices whose distance from $v$ is at most $r$. Since $P_{n,k}$ is an induced subgraph of $W_{n,k}$, the eigenvalues of $P_{n,k}$ interlace those of $W_{n,k}$ 
(see \cite{BH} for eigenvalue interlacing of graphs) in that
\[
\lambda_i(W_{n,k}) \geq \lambda_i(P_{n,k}) \geq \lambda_{i-k}(W_{n,k}). 
\]
This implies that $W_{n,k}$ and $P_{n,k}$ have almost the same eigenvalue distribution. In particular, the eigenvalues of $P_{n,k}$ will be either equal to $\pm 1$ or lie in the intervals $[k-1, k+1]$ and $[-k-1, -k+1]$, with at most $3k$ exceptions. Of the $3k$ possible exceptions, at most $k$ lie in $[1,k-1]$, another $k$ could lie in $[-1,1]$ and the remaining $k$ in $[-k+1, -1]$. 

When $n$ tends to infinity, the eigenvalue distribution of $P_{n,k}$ converges to the spectral distribution of $Z_k$ \cite{M82, MW89}. Since $Z_k$ is vertex-transitive, its whole spectrum is in the essential spectrum, which consists of the set of eigenvalues of infinity multiplicity and the continuous spectrum. Therefore, 
\[
\sigma(Z_k) = \{-1,1\}\cup [k-1, k+1] \cup [-k-1, -k+1].
\]
The spectrum thus consists of two intervals and two additional points. The spectral measure of each interval is $\frac{1}{2k}$, while the two points $\pm 1$ each have spectral measure $\frac{1}{2} - \frac{1}{2k}$. This example partially answers a question of Daniel Lenz and Matthias Keller (private communication), who asked what kind of locally finite vertex-transitive graphs have spectrum which does not consist of a single interval. The motivation to study questions of this kind comes from the study of the spectrum of regular tessellations of hyperbolic spaces. In fact, Schenker and Aizenman \cite{SchenkerA00} found examples of infinite graphs whose spectrum consists of multiple bands. Their examples are not vertex-transitive. 

\bibliographystyle{plain}

\end{document}